\documentclass[11pt]{amsart}
\usepackage[utf8]{inputenc}
\usepackage[english]{babel}
 
\usepackage{color}

\usepackage{amsmath, amsthm, latexsym, amssymb, amsfonts, epsfig, epsf}

\newenvironment{pf}{\proof[\proofname]}{\endproof}
\theoremstyle{plain}
\newtheorem{Th}{Theorem}[section]
\newtheorem{Cor}[Th]{Corollary}
\newtheorem{Conj}[Th]{Conjecture}
\newtheorem{Prop}[Th]{Proposition}
\newtheorem{Lemma}[Th]{Lemma}
\numberwithin{equation}{section} \theoremstyle{definition}
\newtheorem{Rem}[Th]{Remark}

\newtheorem{Def}[Th]{Definition}

\newcommand{\cal}[1]{\mathcal{#1}}

\newcommand{\C}{\mathbb C}
\newcommand{\mS}{\mathbb S}

\newcommand{\R}{\mathbb R}

\newcommand{\D}{\Delta}

\newcommand{\e}{\varepsilon}

\newcommand{\la}{\langle}
\newcommand{\ra}{\rangle}

\newcommand{\supp}{\operatorname{supp}}

\newcommand{\conv}{\operatorname{conv}}

\newcommand{\rs}[1]{Section~\ref{S:#1}}
\newcommand{\rl}[1]{Lemma~\ref{L:#1}}
\newcommand{\rp}[1]{Proposition~\ref{P:#1}}

\newcommand{\re}[1]{(\ref{e:#1})}
\newcommand{\rcj}[1]{Conjecture~\ref{Cj:#1}}

\newcommand{\rt}[1] {Theorem~\ref{T:#1}}

\begin{document}

\title[Wulff shapes and a characterization of simplices]{Wulff shapes and a characterization of simplices via a Bezout type inequality}

\author{Christos Saroglou}
\address[Christos Saroglou]{Department of Mathematical Sciences\\ Kent State University\\ Kent, OH USA}
\email{csaroglo@math.kent.edu}
\author{Ivan Soprunov}
\address[Ivan Soprunov]{Department of Mathematics\\ Cleveland State University\\ Cleveland, OH USA}
\email{i.soprunov@csuohio.edu}
\author{Artem Zvavitch}
\address[Artem Zvavitch]{Department of Mathematical Sciences\\ Kent State University\\ Kent, OH USA}
\email{zvavitch@math.kent.edu}
\thanks{The third author is supported in part by U.S. National Science Foundation Grant DMS-1600753, this material is based upon work supported by the US National Science Foundation under Grant DMS-1440140 while the third author was in residence at the Mathematical Sciences Research Institute in Berkeley, California, during the Fall 2017 semester.}
\date{}

\maketitle

\begin{abstract} Inspired by a fundamental theorem of Bernstein, Kushnirenko, and Khovanskii we study the following Bezout type inequality for mixed volumes
$$
V(L_1,\dots,L_{n})V_n(K)\leq V(L_1,K[{n-1}])V(L_2,\dots, L_{n},K).
$$
We show that the above inequality characterizes simplices, i.e. if $K$ is a convex body satisfying the inequality for all convex bodies $L_1, \dots, L_n \subset \R^n$, then $K$ must be an $n$-dimensional simplex.
The main idea of the proof is to study perturbations given by Wulff shapes. In particular, we prove a new theorem on differentiability of the support 
function of the Wulff shape, which is  of independent interest. 

In addition, we study the Bezout inequality for mixed volumes introduced in \cite{SZ}. We {introduce the} class of weakly decomposable convex bodies which is strictly larger than the set of all polytopes that are non-simplices. We show that the Bezout inequality in \cite{SZ} characterizes weakly indecomposable convex bodies.
\end{abstract}

\section{Introduction}

One of the fundamental results in the theory of Newton polytopes is the Bernstein-Kushnirenko-Khovanskii theorem which expresses algebraic-geometric information such as the degree of an algebraic variety and intersection numbers in terms of convex-geometric invariants such as volumes and mixed volumes \cite{Be, Kho, Kush}. There are many consequences and applications of this beautiful relation in both algebraic and convex geometry (see \cite{KaKh} and references therein). Among them is a recent geometric inequality, called the Bezout inequality for mixed volumes, first introduced in \cite{SZ}. The name comes from interpreting the classical Bezout inequality in algebraic geometry in terms of the mixed volumes with the help of the Bernstein-Kushnirenko-Khovanskii result. In the present paper we study another inequality of this type: For any convex bodies $L_1,\dots, L_n$ in $\R^n$ and any $n$-simplex $K$ one has
\begin{equation}\label{e:B}
V(L_1,\dots,L_{n})V_n(K)\leq V(L_1,K[{n-1}])V(L_2,\dots, L_{n},K).
\end{equation}
Here $V_n(K)$ denotes the $n$-dimensional Euclidean volume and $V(K_1,\dots, K_n)$ denotes the $n$-dimensional mixed volume of bodies 
$K_1,\dots, K_n$. We write $K[{m}]$ to indicate that the body $K$ is repeated $m$ times in the expression for the mixed volume. 
 
The inequality \re{B} has a direct proof based on the monotonicity of the mixed volume presented in \rs{proof}. We also give an algebraic-geometric interpretation of \re{B} at the end of this section. 
Our main result is that \re{B}, in fact, characterizes $n$-simplices in the class of all convex bodies in $\R^n$. In other words, the only convex body 
$K\subset\R^n$ which satisfies \re{B} for all convex bodies $L_1,\dots, L_n$ is the 
$n$-simplex. It turns out that a slightly stronger statement is true.

\begin{Th}\label{T:n-1} Let $K$ be a convex body in $\R^n$. Then $K$ satisfies
\begin{equation}\label{e:main}
V(L_1,\dots,L_{n-1},K)V_n(K)\leq V(L_1,K[{n-1}])V(L_2,\dots, L_{n-1},K,K)
\end{equation}
for all convex bodies $L_1,\dots,L_{n-1}$ in $\R^n$ if and only if $K$ is an $n$-simplex.
\end{Th}

The question of characterization of  $n$-simplices has played an essential role in a number of major results of modern convex geometry, including Ball's {reverse} isoperimetric inequality \cite{Ba} and Zhang's projection inequality \cite{Zh}.  Also simplices are conjecturally the extreme bodies for several open problems. {For example, a {variant} of Bourgain's slicing conjecture \cite{Bo1, Bo2} characterizes simplices as convex bodies with a maximal isotropic constant. The Mahler conjecture states that simplices have the  minimal volume product among all convex bodies. We refer the reader to \cite{Ma1, Ma2, RZ, Sch1} and \cite{Kl} for connections between the slicing problem, the Mahler conjecture, and characterizations of simplices. To the best of the authors' knowledge, (\ref{e:B}) and (\ref{e:main}) are the first characterizations of simplices via an inequality rather than an equality, as 
presented in the above mentioned classical results and open problems.}

Let us briefly describe the method of the proof of \rt{n-1}. The ``if" part follows from \re{B}.
For the ``only if" part, if $K$ is a polytope, the claim follows from the result of \cite{SSZ}, as we show in \rs{polytope}.
In the remaining case, we observe that the problem can actually be seen as an extremal (variational) problem. Indeed, assume that $K$ satisfies (\ref{e:main}) for all convex bodies $L_1,\dots,L_{n-1} \subset \R^n$. Let $K_t$ be a convex perturbation of $K$ where $t\in(-\delta,\delta)$, for some small $\delta>0$, and $K_0=K$. Then the function
$$F(t):=V(K_t,M,K[n-2])V_n(K)-V(K_t,K[{n-1}])V(M,K[{n-1}])$$is non-positive in $(-\delta,\delta)$ and equals 0 at $t=0$, where $M$ is any compact convex set, independent of $t$.
To arrive at a contradiction, one needs to construct $K_t$ and $M$ such that the right derivative of $F$ at $t=0$ (if it exists) is strictly positive. 
The main idea 
is to use a more general construction for the perturbation body $K_t$ than the one used in \cite{SSZ}. 
Namely, given a continuous  function
$f:\mS^{n-1} \to \R$, define $K_t$ as
$$
K_t=\bigcap\limits_{u \in \mS^{n-1}} \left\{ x\in \R^n:  \la x, u\ra \le h_K(u)+tf(u) \right\},
$$
where $h_K$ is the support function of $K$. 

The above construction  goes back to Aleksandrov's 1938 paper \cite{Al1} (see also \cite{Al2, Sch1}) and  utilizes the notion of the Wulff shape.  This notion  was  introduced by G. Wulff \cite{Wu} who studied the equilibrium shape of a droplet or a crystal of fixed volume.  
In the recent years Wulff shapes have been widely used in the study of extremal problems of convex bodies.
We refer to \cite{BLYZ, BLYZ2, MK, HLYZ, HuLYZ, L, Sa, SW} for recent results or \cite{Ga,  Sch1, Se} for expository work on the subject.
A core lemma of Aleksandrov \cite[Lemma 7.5.3]{Sch1} (see also (\ref{e:Aleksandrov-2}) below) provides a differentiability property of the volume of $K_t$.  In our main result in \rs{Wulff} we prove a differentiability property of  $h_{K_t}$ (\rt{dif}): $$
\frac{d h_{K_t}(u)}{dt} \big|_{t=0} = f(u) 
$$
$S_K$-almost everywhere on $\mS^{n-1}$, where $S_K$ is the surface area measure of $K$. 
This can be thought of as a local version of Aleksandrov's lemma. 

There is a close relationship between \re{B} and the Bezout inequality for mixed volumes  studied in \cite{SZ, SSZ}. Recall that it
says that 
\begin{equation}\label{e:Bezout r=2}
V(L_1,L_2,K[{n-2}])V_n(K)\leq V(L_1,K[{n-1}])V(L_2,K[{n-1}]).
\end{equation}
A more general version of this inequality is discussed in \rs{weakly}.
Observe that \re{B} implies \re{Bezout r=2}. It was conjectured in \cite{SZ} that \re{Bezout r=2} characterizes simplices in the class of all convex bodies. This conjecture was confirmed in \cite{SSZ}  for the class of convex $n$-polytopes. Although the general case remains open, the results of this paper provide additional information regarding the conjecture, which we collect in \rs{weakly}.
For example,  \rt{n-1} gives an affirmative answer to the conjecture when $n=3$.
Furthermore, in \rs{weakly} we introduce the notion of weakly decomposable bodies generalizing the classical notion of decomposable bodies. 
The class of weakly decomposable bodies is quite large; we show that it is strictly larger than the set of all polytopes that are non-simplices. In fact, we currently do not have examples of convex bodies which are not weakly decomposable, besides simplices. We show that if $K$ satisfies \re{Bezout r=2} for any $L_1,L_2$ in $\R^n$ then $K$ cannot be weakly decomposable, see \rt{weakly}. 

\rt{n-1} shows that the inequality \re{B} may not hold when $K$ is not an $n$-simplex. Thus, one can ask if \re{B} can be relaxed so that it holds for arbitrary $L_1,\dots, L_n$ and $K$ in $\R^n$. In \rs{Iso} we prove that
$$
V(L_1,\dots,L_{n})V_n(K)\leq n V(L_1,K[{n-1}])V(L_2,\dots, L_{n},K)
$$
for all convex sets $L_1,\dots, L_n$ and $K$ in $\R^n$ and show that this inequality cannot be improved. In addition, we discuss isomorphic versions of \re{main} and \re{Bezout r=2}.

We now turn to an algebraic-geometric interpretation of \re{B}. Recall that the Newton polytope of a polynomial $f\in\C[t_1,\dots,t_n]$ is the convex hull in $\R^n$ of the 
exponent vectors appearing in $f$. Fix $1\leq r\leq n$. A polynomial system $f_1(t)=\dots=f_r(t)=0$ for $t\in(\C\setminus\{0\})^n$ with generic coefficients and fixed Newton polytopes $P_1,\dots, P_r$ defines an $(n-r)$-dimensional algebraic set $X$ in $(\C\setminus\{0\})^n$. 
We may assume that each $P_i$ intersects every coordinate hyperplane $x_i=0$, as translations of the Newton polytopes do not 
change the set $X$. By definition, the degree $\deg X$ of $X$ is the number of intersection points of $X$ with a generic affine subspace of dimension~$r$.  According to the Bernstein-Kushnirenko-Khovanskii theorem, this number can be computed as 
$$\deg X=n!V(P_1,\dots,P_r,\D[n-r]),$$
where $\D=\conv\{0,e_1,\dots, e_n\}$ is the standard $n$-simplex (the Newton polytope of a linear polynomial). Here $e_1,\dots, e_n$ denote the standard basis vectors in $\R^n$.  In particular, when $r=1$, $X$ is a hypersurface and $\deg X$ coincides with the degree of the polynomial defining $X$. When $r=n$, $X$ consists of isolated points whose number equals $\deg X$.  

Now the algebraic-geometric interpretation of \re{B} is as follows. Let $X$ be the hypersurface defined by $f_1(t)=0$ and $Y$ be defined by a polynomial system $f_2(t)=\dots=f_n(t)=0$  in $(\C\setminus\{0\})^n$. As before we assume that the Newton polytopes $P_1,\dots, P_n$ are fixed and the coefficients of the $f_i$ are generic. Then $\deg X=n!V(P_1,\D[n-1])$, $\deg Y=n!V(P_2,\dots, P_n,\D)$, and 
$\deg(X\cap Y)=n!V(P_1,\dots,P_n)$. Also note that $V_n(\D)=1/n!$. The inequality \re{B} then turns into a classical Bezout inequality 
$$\deg(X\cap Y)\leq \deg X\deg Y,$$
see, for instance, Proposition 8.4 of \cite{F} and examples therein.


\section{Preliminaries}\label{sec:pr}

\subsection{Basic definitions} In this section we introduce notation and collect basic facts from classical theory of convex bodies that we use in the paper. As a general reference on the theory we use R. Schneider's book ``Convex bodies: the Brunn-Minkowski theory" \cite{Sch1}.

We write $\la x, y\ra$ for the inner product of $x$ and $y$ in $\R^n$.  Next, $\mS^{n-1}$
denotes the $(n-1)$-dimensional unit sphere in $\R^n$ and $B(x,\delta)$ denotes the closed Euclidean ball of radius
$\delta>0$ centered at $x\in\R^n$. A {\it spherical cap} $U(u,r)$ of radius $r>0$ centered at $u\in\mS^{n-1}$ is the intersection $\mS^{n-1}\cap B(u,r)$. 

For $A\subset\R^n$ its {\it dimension} $\dim A$ is the dimension of the smallest affine subspace of $\R^n$ containing $A$.
A {\it convex body} is a convex compact set with non-empty interior. Note that convex bodies in $\R^n$ are $n$-dimensional.
A {\it (convex) polytope} is the convex hull of a finite set of points. An $n$-dimensional polytope
is called an {\it $n$-polytope} for short. An {\it $n$-simplex} is the convex hull of $n+1$ affinely independent points in $\R^n$. 

For a convex body $K$ the function $h_K:\mS^{n-1}\to\R$, 
$h_K(u)=\max\{\la x, u\ra\ |\ x \in K\}$ is the {\it support function} of $K$. For every convex body $K$ we write $K^u$ to denote the face corresponding to {an outer normal} vector $u \in \mS^{n-1}$, i.e.
$$
K^u=\{ x \in K: \la x, u\ra = h_K(u)\}.
$$
If $\Omega$ is a subset of $\mS^{n-1}$ define the {\it inverse spherical image} $\tau(K,\Omega)$ of $\Omega$ with respect to $K$ by
$$\tau(K,\Omega)=\big\{x\in\partial K:\exists u\in\Omega,\textnormal{ such that }\langle x,u\rangle=h_K(u)\big\}.$$
The {\it surface area measure} $S_K(\cdot)$ of $K$ is a measure on $\mS^{n-1}$ such that
$$S_K(\Omega)={\cal{H}}^{n-1}\big(\tau(K,\Omega)\big),$$
for any Borel  $\Omega \subset \mS^{n-1}$. Here ${\cal{H}}^{n-1}(\cdot)$ stands for the $(n-1)$-dimensional Hausdorff measure.

Let $V_n(K)$  be the Euclidean volume of $K \subset \R^n$.  We will often use the classical formula connecting the volume of a convex body, the support function, and the surface area measure:
\begin{equation}\label{eq:vol}
V_n(K)=\frac{1}{n}\int_{\mS^{n-1}} h_K(u) dS_K(u).
\end{equation}

The Minkowski sum of two sets $K, L \subset \R^n$ is defined as $K+L=\{x+y: x \in K \mbox{  and  } y \in L\}$.  A classical theorem of Minkowski says that if $K_1, K_2, \dots, K_n$ are convex compact sets in $\R^n$ and $\lambda_1, \dots, \lambda_n \ge 0,$ then
$V_n(\lambda_1 K_1 +\lambda_2 K_2+\dots + \lambda_nK_n)$ is a homogeneous polynomial in $\lambda_1,\dots, \lambda_n$ 
of degree $n$.
The coefficient of $\lambda_1\cdots \lambda_n$ is called the {\it mixed volume }of $K_{1}, \dots, K_{n}$ and is denoted by
$V(K_{1}, \dots, K_{n})$. We will also write $V(K_1[m_1],\dots, K_r[m_r])$ for the
mixed volume of  $K_1,\dots, K_r$ where each $K_i$ is repeated $m_i$ times and
$m_1+\dots+m_r=n$. 
We summarize the main properties of the mixed volume below:
\begin{itemize}
\item $V(K,\dots, K)=V(K[n])=V_n(K)$;
\item the mixed volume is symmetric in $K_1,\dots, K_n$;
\item the mixed volume is multilinear: For any $\lambda, \mu \ge 0$ 
$$V(\lambda K + \mu L, K_2,  \dots, K_n)=\lambda V(K,K_2,  \dots, K_n)+ \mu V(L,K_2, \dots, K_n);$$
\item the mixed volume is translation invariant: For any $a \in \R^n$
$$V(K_1+a,K_2,  \dots K_n)= V(K_1,K_2, \dots, K_n);$$ 
\item the mixed volume is monotone: If $ K \subset L$ then $$V(K,K_2, \dots, K_n) \le V(L,K_2, \dots, K_n).$$ 
\end{itemize}

We will need the following classical inequalities for the mixed volumes: 
 The Minkowski First inequality
 \begin{equation}\label{eq:mink}
 V(K, L[n-1])\ge V_n(K)^{1/n}V_n(L)^{(n-1)/n},
\end{equation}
and the Aleksandrov--Fenchel inequality
\begin{equation}\label{eq:af}
V(K_1,K_2, K_3, \dots, K_n) \ge \sqrt{V(K_1,K_1, K_3, \dots, K_n)V(K_2,K_2, K_3, \dots, K_n)}.
\end{equation}


Let $S(K_1,\dots, K_{n-1},\cdot)$ be the {\it mixed area measure} for bodies $K_1,\dots, K_{n-1}$ defined by
\begin{equation}\label{eq:mixarea}
V(L,K_1,\dots, K_{n-1})=\frac{1}{n}\int_{\mS^{n-1}}h_L(u)dS(K_1,\dots, K_{n-1},u)
\end{equation}
for any compact convex set $L$. In particular,
\begin{equation}\label{eq:firstmix}
V(L,K[n-1])=\frac{1}{n}\int_{\mS^{n-1}}h_L(u)dS_K(u)
\end{equation}
and
$$
 S(K,\dots,K,\cdot)=S(K[n-1],\cdot)=S_K(\cdot).
 $$
The identity (\ref{eq:mixarea}) implies that the invariance properties of mixed volumes are inherited by mixed area measures. More specifically, $S(K_1,\dots,K_{n-1},\cdot)$ is translation invariant, symmetric, and multilinear with respect to $K_1,\dots,K_{n-1}$.
We also note that if the $K_i$ are polytopes the mixed area measure $S(K_1,\dots, K_{n-1},\cdot)$ has finite support and
for every $u\in\mS^{n-1}$ we have
\begin{equation}\label{e:mixed-area}
S(K_1,\dots, K_{n-1},\{u\})=V(K_1^u,\dots, K_{n-1}^u),
\end{equation}
where $V(K_1^u,\dots, K_{n-1}^u)$ is the $(n-1)$-dimensional mixed volume of the faces $K_i^u$ translated 
the the subspace orthogonal to $u$,
see \cite[Sec 5.1]{Sch1}. 
We will need a slightly more general statement. 
\begin{Lemma}\label{L:mvol} Let $K_1,\dots, K_{n-1}$ be convex bodies in $\R^n$. Then
$$S(K_1,\dots,K_{n-1},\{u\})=V(K^u_1,\dots,K^u_{n-1}).$$
\end{Lemma}

\begin{pf} We use an alternative definition of the mixed area measure given in \cite[(5.21), page 281]{Sch1}:
$$S(K_1,\dots, K_{n-1},\cdot)=\frac{1}{(n-1)!}\sum_{k=1}^{n-1}(-1)^{n-1-k}\sum_{i_1<\dots<i_k}S_{K_{i_1}+\dots+K_{i_k}}(\cdot).$$
By the definition of the surface area measure
$$S(K_{i_1}+\dots+K_{i_k},\{u\})=V_{n-1}((K_{i_1}+\dots+K_{i_k})^u),$$
where as before $(K_{i_1}+\dots+K_{i_k})^u$ denotes the face of $K_{i_1}+\dots+K_{i_k}$ corresponding to $u$. Note
that $(K_{i_1}+\dots+K_{i_k})^u=K^u_{i_1}+\dots+K^u_{i_k}$.
Therefore,
$$S(K_1,\dots, K_{n-1},\{u\})=\frac{1}{(n-1)!}\sum_{k=1}^{n-1}(-1)^{n-1-k}\sum_{i_1<\dots<i_k}V_{n-1}(K^u_{i_1}+\dots+K^u_{i_k}).$$
The right hand side of the above equation is the well-known formula for the mixed volume $V(K^u_1,\dots,K^u_{n-1})$, see \cite[Lemma 5.1.4]{Sch1}.
\end{pf}

Finally, we will need the following formula relating the mixed volumes and projections, which is a particular case of \cite[Theorem 5.3.1]{Sch1}.
For $u\in\mS^{n-1}$ the orthogonal projection of a set $A\subset\R^n$ onto the subspace  $u^{\bot}$ orthogonal to $u$ is denoted by $A|u^{\bot}$.
Let $I=[0,u]$ be a unit segment for some $u\in\mS^{n-1}$ and $K_2,\dots, K_n$ convex bodies in $\R^n$. Then 
\begin{equation}\label{e:proj}
V(I,K_2,\dots, K_n)=\frac{1}{n}V(K_2|u^{\bot},\dots,K_n|u^{\bot}),
\end{equation}
where in the right-hand side is the $(n-1)$-dimensional mixed volume of convex bodies in the orthogonal subspace $u^\bot$.


\section{{ Differentiation of Wulff Shapes.}}\label{S:Wulff}


The notion of the Wulff shape is the main tool in constructing the perturbation $K_t$. We refer to Section 7.5 of \cite{Sch1} for details.
We start with the definition of the Wulff shape and several properties that we need later.

\begin{Def}
Consider a non-negative function $g: \mS^{n-1} \to \R$. The {\it Wulff shape} $W(g)$ of $g$ is defined as
$$
W(g)=\bigcap\limits_{u \in \mS^{n-1}} \left\{ x\in \R^n:  \la x, u\ra \le g(u) \right\}.
$$
\end{Def}
Note that when $g$ is the support function $h_K$ of a convex body $K$ then $W(h_K)=K$.
\begin{Prop}\label{P:wc} Consider non-negative functions $g_1, g_2: \mS^{n-1} \to \R$ and $\lambda \in (0,1)$. Then
$$
W((1-\lambda)g_1 + \lambda g_2)  \supseteq (1-\lambda) W(g_1)+ \lambda W (g_2).
$$
\end{Prop}
\pf  Indeed consider $x \in W(g_1)$ and $y\in W(g_2)$. Then 
$$
\la x, u\ra \le g_1(u)  \mbox{   and  }   \la x, u\ra \le g_2(u)
$$
for all $u \in \mS^{n-1}$. Thus
$$
(1-\lambda)x+\lambda y \in \{ z\in \R^n: \la z, u\ra \le (1-\lambda)g_1(u)+ \lambda g_2(u)\}
$$
for all $u\in \mS^{n-1}$, which completes the proof.
\endpf
Consider a convex body $K\subset \R^n$, containing the origin in its interior. Given a continuous  function 
$f:\mS^{n-1} \to \R$, we define the {\it perturbation} $K_t(f)$ of $K$ as
\begin{equation}\label{e:perturb}
K_t(f)=W(h_K+ tf),
\end{equation}
where $t \in (a, b)$ such that $h_K(u)+tf(u) >0$ for all $u \in \mS^{n-1}$. Note that $K_t(f)$ is also a convex body containing the origin in its interior.
To shorten our notation we write $K_t$ instead of $K_t(f)$ when it is clear what function $f$ is being considered.

The following proposition shows that the support function of the perturbation $K_t$ is concave with respect to $t$.

\begin{Prop}\label{P:wh} The support function $h_{K_t}(u)$ is concave with respect to $t\in (a,b)$ for all $u \in \mS^{n-1}$.
\end{Prop}
\pf Consider $\lambda \in [0,1]$ and $t,s \in (a,b)$. Then, using \rp{wc}, we get
\begin{equation*}
\begin{split}
K_{\lambda t+ (1-\lambda)s} &= W\left(\lambda h_K +(1-\lambda)h_K +(\lambda t + (1-\lambda)s) f\right)\\ &\subseteq \lambda W(h_K+ tf) +(1-\lambda)W(h_K+ sf).
\end{split}
\end{equation*}
\endpf

Aleksandrov's lemma \cite[Lemma 7.5.3]{Sch1} states that if $h:\mS^{n-1}\to(0,\infty)$ and $f:\mS^{n-1}\to\R$ are continuous functions, then 
\begin{equation}\label{e:Aleksandrov}
\frac{d}{dt}W(h+tf)|_{t=0}=\int\limits_{\mS^{n-1}} f(u) d S_{W(h)}(u).
\end{equation}
In particular, when $h=h_K$ one has
\begin{equation*}\label{e:Aleksandrov-2}
\frac{d}{dt} V(K_t)|_{t=0}=\int\limits_{\mS^{n-1}} f(u) d S_K(u).
\end{equation*}
Aleksandrov's lemma (\ref{e:Aleksandrov}) follows from the following fact (see the proof of Lemma 7.5.3 in \cite{Sch1}), which is going to be crucial for the proof of Theorem \ref{T:dif} below.
\begin{Prop}\label{P:der}
Let $h:\mS^{n-1}\to(0,\infty)$ and $f:\mS^{n-1}\to\R$ be continuous functions. Then, 
$$\frac{d}{dt}V(W(h+tf), W(h), \dots, W(h))|_{t=0}=\frac{1}{n} \int\limits_{\mS^{n-1}} f(u) d S_{W(h)}(u).$$
In particular, when $h=h_K$
$$
\frac{d}{dt}V(K_t(f), K, \dots, K)|_{t=0}=\frac{1}{n} \int\limits_{\mS^{n-1}} f(u) d S_K(u).
$$
\end{Prop}
%

\noindent Now we go back to the support function $h_{K_t}(u)$. By concavity established in \rp{wh}, the one-sided derivatives of $h_{K_t}(u)$ exist at every $t$. In particular, for $t=0$ we have the following theorem.

{\begin{Th}\label{T:dif}
$$
\frac{d h_{K_t}(u)}{dt} \big|_{t=0} = f(u) 
$$
$S_K$-almost everywhere on $\mS^{n-1}.$
\end{Th}
\pf
We need to show that
\begin{equation}\label{e:S_K-a.e.}
\lim\limits_{t\to 0^+}  \frac{h_{K_t}(u) - h_K(u) }{t}= f(u)
\quad\text{and}\quad
\lim\limits_{t\to 0^-}  \frac{h_{K_t}(u) - h_K(u)}{t} = f(u),
\end{equation}
$S_K$-almost everywhere on $\mS^{n-1}.$
Indeed, it follows from 
the definition of the Wulff shape that
$$
\lim\limits_{t\to 0^+}  \frac{h_{K_t}(u) - h_K(u)}{t} \le f(u)
$$
for all $u \in \mS^{n-1}$. Assume that there is a set of positive $S_K$-measure where we have a strict inequality, then
$$
\int\limits_{\mS^{n-1}}\lim\limits_{t\to 0^+}  \frac{h_{K_t}(u) - h_K(u)}{t}  dS_K(u) <  \int\limits_{\mS^{n-1}} f(u) dS_K(u).
$$
On the other hand, \rp{der}, together with  the Dominating Convergence theorem gives
$$
\frac{1}{n} \int\limits_{\mS^{n-1}} f(u) dS_K(u) =\frac{1}{n}\lim\limits_{t\to 0+} \int\limits_{\mS^{n-1}}  \frac{h_{K_t}(u) - h_K(u)}{t}  dS_K(u)
$$
$$
= \frac{1}{n} \int\limits_{\mS^{n-1}}\lim\limits_{t\to 0^+}  \frac{h_{K_t}(u) - h_K(u)}{t}  dS_K(u) < \frac{1}{n} \int\limits_{\mS^{n-1}} f(u) dS_K(u).
$$
which gives a contradiction. We use similar logic to prove the second part of (\ref{e:S_K-a.e.}).
\endpf}

\begin{Rem}\label{R:sharp}
In view of Theorem \ref{T:dif}, one might ask if  $\frac{d h_{K_t}(u)}{dt} \big|_{t=0}=f(u)$ holds for all  $u \in \mS^{n-1}$. We note that the answer to this question is negative, in general. Indeed, let $K=\conv\{0, e_1, e_2\}$ be the standard triangle in $\mathbb{R}^2$. 
Let $f:\mS^1 \to [0,1]$ be a continuous function such that $f(e_2)=1$ and $f$ is zero outside of a small neighborhood of $e_2$. Then $K_t=K,$ for $t \ge 0$ and thus $\lim\limits_{t \to 0^+} \frac{h_{K_t} -h_K}{t}  =0$. Moreover, if $-1<t<0$  then
$h_{K_t}(e_1) \le 1+t $ and thus $\lim\limits_{t \to 0^-} \frac{h_{K_t}(e_2) -h_K(e_2)}{t}  \ge 1$. Therefore, the derivative $\frac{d h_{K_t}(e_2)}{dt} \big|_{t=0}$ does not exist. 
\end{Rem}


\section{Proof of \rt{n-1}}\label{S:proof}

We begin with the ``if" part of the statement. Let $K$ be an $n$-simplex and $L_1,\dots,L_n$ convex bodies in $\R^n$. Then
\begin{equation}\label{e:BB}
V(L_1,\dots,L_{n})V_n(K)\leq V(L_1,K[{n-1}])V(L_2,\dots, L_{n},K).
\end{equation}

Indeed, since this inequality is invariant under dilations and translations of $L_1$ we may assume that 
$L_1$ is contained in $K$ and intersects every facet of $K$. Then, by (\ref{eq:firstmix}), we have 
\begin{equation}\label{e:BB proof}
V(L_1,K[{n-1}])=\frac{1}{n}\int\limits_{\mS^{n-1}}h_{L_1}(u)d S_K(u)=\frac{1}{n}\sum_{i=0}^n h_{L_1}(u_i)V_{n-1}(K^{u_i}),
\end{equation}
where $\{u_0,\dots,u_n\}$ is the set of the outer unit normals to the facets of $K$. But $h_{L_1}(u_i)=h_K(u_i)$ since $L_1\subseteq K$ and it intersects each facet of $K$. Therefore, the last expression in \re{BB proof} equals $V_n(K)$. Now \re{BB proof} follows from the monotonicity 
property of the mixed volume, 
$$V(L_1,\dots,L_{n})\leq V(L_2,\dots, L_{n},K).$$
Finally, we note that \re{main} is a particular case of \re{BB}.

To prove the ``only if'' part 
we treat the following three cases separately: (1) $K$ is a polytope, 
(2) $K$ has an infinite number of facets, and (3) $K$ has a finite number of facets but is not a polytope. 
Then \rt{n-1} follows from Theorems \ref{T:polytope}, \ref{T:case2}, and~\ref{T:case1} below.


 \subsection{The case of polytopes.}\label{S:polytope} Suppose $K$ is  a polytope which satisfies
 \begin{equation}\label{e:main-1}
V(L_1,\dots,L_{n-1},K)V_n(K)\leq V(L_1,K[{n-1}])V(L_2,\dots, L_{n-1},K,K)
\end{equation}
 for any bodies $L_1,\dots, L_{n-1}$. In particular, when $L_3=\dots=L_{n-1}=K$ we obtain
  \begin{equation}\label{e:polytope}
V(L_1,L_2,K[{n-2}])V_n(K)\leq V(L_1,K[{n-1}])V(L_2,K[n-1]).
\end{equation}
Now the fact that $K$ is an $n$-simplex follows from \cite[Theorem~3.4]{SSZ} which we formulate next.

\begin{Th}\label{T:polytope}
Let $K$ be an $n$-polytope in $\R^n$. Suppose  \re{polytope}
holds for all convex bodies $L_1$ and $L_2$ in $\R^n$. Then $K$ is a simplex.
\end{Th}
{Theorem \ref{T:polytope} follows also from the results about weakly decomposable convex bodies in Section \ref{S:weakly} below.}


\subsection{$K$ has an infinite number of facets.} The following result implies \rt{n-1} in this case.

\begin{Th}\label{T:case2}
Let $K$ be a convex body in $\R^n$ with infinitely many facets. Then there exist 
convex bodies $L_1$ and $L_2$ in $\R^n$ such that the inequality 
\begin{equation}\label{e:case2}
V(L_1,L_2,K[n-2])V_n(K)\leq V(L_1,K[{n-1}])V(L_2,K[{n-1}])
\end{equation}
is false.
\end{Th}

\begin{pf} 
Since $K$ has infinitely many facets, {for any $\e>0$, there exists a facet $K^u\subset K$ such that $V_{n-1}(K^u) \le \e$, otherwise the surface area of $K$ would be infinite which contradicts the definition of a convex body. Set $\e_0=V_{n-1}(K^u)\leq\e$.} 
Since $K$ has non-empty interior it contains a ball of radius $\delta>0$, independent of $\e_0$. 
Cut this ball in half by a hyperplane orthogonal to $u$ and let $M$ be the half whose facet has outer normal $u$. 
Furthermore, let $K_t=K_t(f)$ be a perturbation as in \re{perturb} for some {continuous } $f:\mS^{n-1}\to[0,1]$ and $t\in(-a,a)$ {which we choose  below}. 

Now assume \re{case2} holds for any $L_1,L_2$. Put $L_1=K_t$ and $L_2=M$ and consider the function
$$F(t)=V(K_t,M,K[n-2])V_n(K)-V(K_t,K[{n-1}])V(M,K[{n-1}]).$$
Then $F(0)=0$ and $F(t)\leq 0$ on $(-a,a)$. Our goal is to show that for small enough $\e$ the right-sided derivative of $F(t)$ at $t=0$ is positive,
which provides a contradiction.
We deal with each summand in the definition of $F(t)$ separately. Using (\ref{eq:mixarea}) we get 
\begin{equation*}
\begin{split}
\lim\limits_{t\to 0+}\frac{1}{t} \big(V(K_t, M,K[n-2])&-V(K, M, K[n-2])\big) \\
&=\frac{1}{n} \lim\limits_{t\to 0+}\frac{1}{t} \int\limits_{\mS^{n-1}} \left(h_{K_t}(x)-h_K(x)\right) d S(M,K[n-2],x)\\
&\ge\frac{1}{n} \lim\limits_{t\to 0+}\frac{1}{t} \int\limits_{\{u\}}\left(h_{K_t}(x)-h_K(x)\right) d S(M,K[n-2],x) \\
&= \frac{1}{n} f(u)V(M^u,K^u[n-2]),
\end{split}
\end{equation*}
where the last equality follows from and (\ref{e:mixed-area}) and \rt{dif} (note that since $S_K(\{u\})=\e_0>0$, the derivative of $h_{K_t}(u)$ at $t=0$ exists).
Setting $f(u)=1$ and using the Minkowski inequality (\ref{eq:mink})  for the mixed volume $V(M^u,K^u[n-2])$ we get
$$
\lim\limits_{t\to 0+}\frac{1}{t} \big(V(K_t, M,K[n-2])-V(K, M, K[n-2])\big) \ge \frac{1}{n}V_{n-1}(M^u)^{\frac{1}{n-1}}V_{n-1}(K^u)^{\frac{n-2}{n-1}}.
$$
Recall that $V_{n-1}(K^u)=\e_0$ and $V_{n-1}(M^u)$ is positive, independent of $\e_0$. Thus, we can write
\begin{equation}\label{e:LHS}
\lim\limits_{t\to 0+}\frac{1}{t} \big(V(K_t, M,K[n-2])-V(K, M, K[n-2])\big) \ge C_1\e_0^{\frac{n-2}{n-1}},
\end{equation}
for some $C_1>0$, independent of $\e_0$.

Now we turn to the second summand in $F(t)$. This is straightforward from \rp{der}:
\begin{equation}\label{e:RHS}
\lim\limits_{t\to 0^+}\frac{1}{t}\left(V(K_{t},K[n-1])-V_n(K)\right)=\frac{d}{dt}V(K_{t},K[n-1])|_{t=0}=\frac{1}{n}\int\limits_{\mS^{n-1}} f(x) d S_K(x).
\end{equation}
We can choose $f$ to be positive on the facet $K^u$ and zero outside. Then the integral above equals $C_2\e_0$ for some $C_2>0$, independent of  $\e_0$. Combining \re{LHS} and \re{RHS} we obtain
$$
\lim\limits_{t\to 0^+} \frac{F(t)}{t}\geq C_1V_n(K)\e_0^{\frac{n-2}{n-1}}-C_2V(M,K[n-1])\e_0>0,
$$
for small enough $\e_0$ and, consequently, for small enough $\e$. This completes the proof of the theorem.
\end{pf}



 \subsection{$K$ has finitely many facets and is not a polytope.}  
Although the methods in this case are similar to the ones used in the previous case, the proof is more involved.  

Let  $\supp(S_K)$ be the support of the measure $S_K$, i.e.  the largest (closed) subset of $\mS^{n-1}$ for which every open neighborhood of every point of the set has positive measure.  We start with the following observation.

\begin{Prop}\label{P:bad points}
Let $K$ be a convex body in $\R^n$ which satisfies 
  \begin{equation}\label{e:polytope-2}
V(L_1,L_2,K[{n-2}])V_n(K)\leq V(L_1,K[{n-1}])V(L_2,K[n-1])
\end{equation}
for all $L_1,L_2$. Then $\dim K^u>0$ for every $u\in\supp(S_K)$.
\end{Prop}

\begin{pf} The proof is similar to the one of \cite[Theorem 4.1]{SSZ}. Assume there exists $u\in\supp(S_K)$ such that $K^u=\{y\}$, for some $y\in \partial K$. 
For $\varepsilon>0$ define the truncated body
$$K_\varepsilon=\{x\in K : \langle x, u\rangle \leq h_K(u)-\varepsilon\}.$$
Since $K^u=\{y\}$, it follows that there exist 
$v\in \mS^{n-1}$ and $\varepsilon>0$ such that 
the projections of $K$ and $K_\varepsilon$ are the same: 
\begin{equation}\label{e:last1}
K_\varepsilon|v^{\bot}=K|v^{\bot}. 
\end{equation}
A proof of this statement can be found in the proof of \cite[Theorem 4.1]{SSZ}.

Let $L_1=[0,v]$ and $L_2=K_\varepsilon$. Then, from \re{proj} and (\ref{e:last1}), we get 
$$V(L_1,L_2,K[n-2])=\frac{1}{n}V(K_\varepsilon|v^{\bot},K|v^{\bot}[n-2])=\frac{1}{n}V(K|v^{\bot})=V(L_1,K[n-1]).$$
Therefore, it suffices to show that $V_n(K)>V(L_2,K[n-1])$; this will contradict our assumption (\ref{e:polytope-2}). Note
that $h_{K_\varepsilon}\leq h_K$ and 
$$
V(K)-V(L_2,K[n-1])=\frac{1}{n}\int_{\mS^{n-1}}(h_K(v)-h_{K_\varepsilon}(v))dS_K(v),
$$
by  (\ref{eq:vol}) and (\ref{eq:firstmix}). Therefore
it is enough to show that there exists a Borel set $\Omega\subseteq \mS^{n-1}$, such that $S_K(\Omega)>0$ and $h_K(v)>h_{K\varepsilon}(v)$, for all $v\in\Omega$. To this end we claim the following.

\textit{Claim.} There exists an open set $U\subseteq \mS^{n-1}$, such that $U$ contains $u$ and $\tau(K,U)\subseteq K\setminus K_\varepsilon$. To see this, assume that this is not the case. Then, there is a sequence of open spherical caps $U_N=U(u,r_N)$ of radii $r_N\to 0$, as $N\to\infty$, such that $\tau(K,U_N)$ contains a point $x_N\in K_\varepsilon$. This implies that there exists $u_N\in U_N$ such that $u_N\in K^{u_N}$. Since $u_N\to u$, by compactness there exists $x\in \partial K\cap K_\varepsilon$, such that $x\in K^u$. This contradicts our assumption
$K^u=\{y\}$, since $x\neq y$ as $y\not\in K_\varepsilon$.

An immediate reformulation of the previous claim states that there exists an open set $U$ in $\mS^{n-1}$ containing $u$, such that $K^v\subseteq K\setminus K_\varepsilon$, for all $v\in U$. In other words, $h_K(v)>h_{K_\varepsilon}(v)$, for all $v\in U$. However, since $u\in \supp(S_K)$, we have $S_K(U)>0$, which completes the proof.
\end{pf}

Suppose $K$ satisfies the conditions of \rp{bad points}. Then we can choose an open subset 
 $U  \subseteq \mS^{n-1}$ satisfying the following three conditions:
\begin{enumerate}
\item[(i)] $V=U \cap \supp(S_K) \neq\varnothing$ (and, hence, $S_K(V)>0$),
\item[(ii)] $V$ contains no normals to facets of $K$,
\item[(iii)] there exists $\ell>0$ such that for any $u\in V$ the face $K^u$ contains a segment of length at least $\ell$.
\end{enumerate}
Indeed, we can ensure (i) and (ii) are satisfied since we can exclude the finite set of unit normals to the facets of $K$ from the support of $S_K$.
For part (iii), choose $V'$ satisfying (i) and (ii) and consider the sequence of subsets
$$V_i=\{u\in V':  K^u \text{ contains a segment of length at least } {1}/{i}\},\ \ i\in\mathbb{N}.$$
By \rp{bad points}, $V_i$ is nonempty for $i$ large enough and  $\cup_{i\geq 1}V_i=V'$.
Therefore, there exists $m\geq 1$ such that $V_m$ has positive $S_K$-measure. Put $\ell=1/m$ and $V=V_m$.

Our next step is to construct a sequence of perturbations $K_t(f_N)$ for certain functions $f_N$. We begin by
constructing a nested sequence $U_N\subset U$ of spherical caps centered at points of $V$.

Let $U_1=U(u_0,\rho_1)\subset U$ for some $u_0\in V$ and $\rho_1>0$. Consider a continuous function $f_1:\mS^{n-1}\to [0,1]$ such that $f_1$
is zero on $\mS^{n-1}\setminus U_1$, positive on $U_1$, and takes values at least~$\frac{1}{2}$ on ${U_1'}:=U(u_0,\frac{\rho_1}{2})$.
By \rt{dif} there exists $u_1\in {U_1'}\cap V$ such that $h_{K_t(f_1)}$ has derivative equal to $f_1(u_1)$.
Define $U_2=U(u_1,\rho_2)\subset U_1$ where $\rho_2\leq\frac{\rho_1}{2}$ and continue the process. We obtain a sequence of spherical caps $U_N=U(u_{N-1},\rho_N)$ with centers at $u_{N-1}\in V$ and radii $\rho_N\to 0$, as $N\to\infty$, and a sequence of functions $f_N:\mS^{n-1}\to [0,1]$ satisfying:
\begin{enumerate}
\item[(a)] $f_N$ is positive on $U_N$, 
\item[(b)] $f_N$ is zero on $\mS^{n-1}\setminus U_N$,
\item[(c)] $\frac{d}{dt}h_{K_t(f_N)}|_{t=0}=f_N(u_{N})\geq\frac{1}{2}$.
\end{enumerate}


We are ready for the main result in the case when $K$ is not a polytope and has {at most} finitely many facets.

\begin{Th}\label{T:case1}
Let $K$ be a convex body in $\R^n$ with finitely many facets which is not a polytope. 
Then there exist convex bodies $L_1,\dots, L_{n-1}$ in $\R^n$ such that the inequality 
\begin{equation}\label{e:case1}
V(L_1,\dots,L_{n-1},K)V_n(K)\leq V(L_1,K[{n-1}])V(L_2,\dots, L_{n-1},K,K)
\end{equation}
is false.
\end{Th}

\begin{pf} Suppose $K$ satisfies \re{case1} for all $L_1,\dots, L_{n-1}$. In particular, it satisfies the conditions of \rp{bad points} and,
hence, there exist a set $V\subset\mS^{n-1}$, a sequence of points $u_N\in V$, and perturbations $K_{t,N}=K_t(f_N)$ as
 constructed above. Choose a ball $B$ of radius $\delta>0$ contained in $K$ and cut it in half by a hyperplane orthogonal to $u_N$.
 Let $M_N$ be the half of $B$ whose facet has outer normal $u_N$. We set $L_1=K_{t,N}$ and $L_i=M_N$ for $2\leq i\leq n-1$.
 Then \re{case1} produces
\begin{equation}\label{e:B2case1}
V(K_{t,N},M_N[n-2],K)V_n(K)\leq V(K_{t,N},K[n-1])V(M_N[n-2],K,K).
\end{equation}
For every $N\geq 1$, define the function
\begin{equation}\label{e:function}
F_N(t)=V(K_{t,N},M_N[n-2],K)V_n(K)- V(K_{t,N},K[n-1])V(M_N[n-2],K,K).
\end{equation}
Note that $F_N(0)=0$ and if \re{B2case1} holds then $F_N(t)\leq 0$ for all $t\in(-a,a)$. 
As in the proof of \rt{case2} we show that for $N$ large enough the right-sided derivative of $F_N(t)$ at $t=0$ is positive,
which gives a contradiction. 

For the first summand in $F(t)$ we have
\begin{equation*}
\begin{split}
V(K_{t,N},M_N[n-2],K)&-V(K,M_N[n-2],K)\\
&=\frac{1}{n}\int\limits_{\mS^{n-1}}\left(h_{K_{t,N}}(x)-h_K(x)\right)dS(M_N[n-2],K,x)\\
&\geq\frac{1}{n}\int\limits_{\{u_N\}}\left(h_{K_{t,N}}(x)-h_K(x)\right)dS(M_N[n-2],K,x)\\
&=\frac{1}{n}\left(h_{K_{t,N}}(u_N)-h_K(u_N)\right)V(M^{u_N}_N[n-2],K^{u_N}),
\end{split}
\end{equation*}
where the last equality follows \rl{mvol}.  We claim that the mixed volume is bounded below by a positive constant $C$ independent of $N$:
$$V(M^{u_N}_N[n-2],K^{u_N})\geq C>0.$$
Indeed, by condition (iii) in the definition of $U$, the face $K^{u_N}$ contains a segment $I$ of length at least $\ell>0$. Also, by construction, the facet $M^{u_N}_N$ is an $(n-1)$-dimensional disc $D_\delta$ of radius $\delta>0$. Thus, 
$$V(M^{u_N}_N[n-2],K^{u_N})\geq V(D_\delta[n-1],I)=\frac{\ell}{n}V_{n-1}(D_\delta)=C>0.$$
Therefore, for the first summand of $F_N(t)$ we have
\begin{equation*}
\begin{split}
\lim\limits_{t\to 0^+}\frac{1}{t}\big(V(K_{t,N},&M_N[n-2],K)-V(K,M_N[n-2],K)\big) \\
&\geq C\lim\limits_{t\to 0^+}\frac{h_{K_{t,N}}(u_N)-h_K(u_N)}{t}=Cf_N(u_N)\geq\frac{C}{2}.
\end{split}
\end{equation*}
The last equality follows from Theorem \ref{T:dif} and the last inequality follows from condition (c) in the definition of $f_N$.

For the second summand of $F_N(t)$ we use \rp{der} as before:
$$
\lim\limits_{t\to 0^+}\frac{1}{t}\left(V(K_{t,N},K[n-1])-V_n(K)\right)=\frac{d}{dt}V(K_{t,N},K[n-1])|_{t=0}=\frac{1}{n}\int\limits_{\mS^{n-1}} f_N(x) d S_K(x).
$$
Bringing the two summands together we obtain
$$
\lim\limits_{t\to 0^+} \frac{F_N(t)}{t}\geq \frac{C}{2}V_n(K)-\frac{1}{n}V(M_N[{n-2}],K,K)\int\limits_{\mS^{n-1}} f_N(x) d S_K(x).
$$
Recall that $f_N$ is bounded and vanishes outside of a spherical cap of radius $\rho_N\to 0$. By choosing $N$ large enough we can ensure that the integral of $f_N$ is smaller than any given number. On the other hand, $V(M_N[{n-2}],K,K)\leq V_{n}(K)$ since $M_N\subset K$, for any $N\geq 1$. Therefore, it is enough to take $N$ large enough so 
$\int\limits_{\mS^{n-1}} f_N(x) d S_K(x)<nC/2$ to obtain
$$
\lim\limits_{t\to 0^+} \frac{F_N(t)}{t}>0.
$$
This completes the proof of the theorem.
\end{pf}


\section{Weakly decomposable convex bodies}\label{S:weakly}

In this section we discuss the impact of the results of \rs{proof} on the conjecture about the Bezout inequality for mixed volumes
formulated in \cite{SZ}. Recall this conjecture:

 \begin{Conj}\label{Cj:BE} Fix an integer $2\leq r\leq n$ and 
 let $K\subset\R^n$ be a convex body satisfying the Bezout inequality 
\begin{equation}\label{e:Bezout}
V(L_1,\dots,L_r,K[{n-r}])V_n(K)^{r-1}\leq \prod_{i=1}^r V(L_i,K[{n-1}])
\end{equation}
for all convex bodies $L_1,\dots, L_r$ in $\R^n$. Then $K$ is an $n$-simplex.
 \end{Conj}

We then introduce a class of weakly decomposable convex bodies which generalizes the classical notion of decomposable bodies.
We show that every polytope which is not a simplex is weakly decomposable and there are many weakly decomposable bodies which are not polytopes. Our main result of this section (\rt{weakly}) asserts that a convex body $K$ satisfying \re{Bezout} for all 
convex bodies $L_1,\dots, L_r$ cannot be weakly decomposable, which improves  the result of \cite[Theorem 3.3]{SZ}.

\subsection{Impact of results of \rs{proof} on \rcj{BE}.}
The special case of the inequality \re{Bezout} with $r=2$ was already considered in the introduction, see \re{Bezout r=2}.
As we mentioned there, \re{B} implies \re{Bezout r=2}. In fact, \re{B} implies \re{Bezout} for any $2\leq r\leq n$. Indeed, \re{Bezout} with $r=n$ can be obtained from \re{B} by successive iterations. Clearly the case of $r=n$ implies all the other cases.

Next, observe that if \rcj{BE} is true for $r=2$ then it is true for any $2\leq r\leq n$. In dimension three the conjecture (with $r=2$)
says that if a convex body $K\subset\R^3$ satisfies 
\begin{equation}\label{e:Bezout n=3}
V(L_1,L_2,K)V_3(K)\leq V(L_1,K,K)V(L_2,K,K)
\end{equation}
for any convex bodies $L_1,L_2$ in $\R^3$ then $K$ is a $3$-simplex. This is precisely the statement of  \rt{n-1} with $n=3$. Therefore, \rcj{BE} is true in the three-dimensional case. In the general case, the results of \rt{polytope}, \rp{bad points}, and \rt{case2} provide with the following new information regarding \rcj{BE}.

\begin{Cor} Let $K\subset\R^n$ for $n\geq 4$ be a convex body which is not a polytope and which satisfies \re{Bezout} for 
all convex bodies $L_1,\dots, L_r$ in $\R^n$. Then $K$ has at most finitely many facets and for every $u$ in the support of the surface area measure
$S_K$, the face $K^u$ is positive dimensional.
\end{Cor} 

\subsection{Weakly decomposable bodies}
Recall that a convex body is called \textit{decomposable} if $K=L+M$, for some compact convex sets $L$, $M$ which are not homothetic to~$K$. It was shown in \cite{SZ} that decomposable convex bodies do not satisfy the Bezout inequality \re{Bezout} for any $2\leq r\leq n$. We generalize the definition of decomposability as follows.
\begin{Def}\label{D:decomposable}
A convex body $K$ in $\mathbb{R}^n$ is called \textit{weakly decomposable} if there exists a convex set $M$, not homothetic to $K$, such that the surface area measure $S_{K+M}$ of $K+M$ is absolutely continuous with respect to the surface area measure $S_K$ of $K$. A convex body which is not weakly decomposable is called  \textit{weakly indecomposable}.
\end{Def}

The following proposition justifies the terminology.

\begin{Prop}\label{P:decomp} 
Every decomposable convex body is weakly decomposable.
\end{Prop}

To prove this we need the following lemma.

\begin{Lemma}\label{l:weakly decomposable}
Let $M$ and $L$ be compact convex sets in $\mathbb{R}^n$. Then
\begin{enumerate}
\item[i)] $\displaystyle{S_{L+M}=\sum_{r=0}^{n-1}\binom{n-1}{r}S(M[r],L[n-r],\cdot)}$,
\item[ii)] $S_{L+M}$ is absolutely continuous with respect to $S_L$ if and only if the mixed area measure $S(M[r],L[n-r],\cdot)$ is absolutely continuous with respect to $S_L$ for all $0\leq r<n$.
\end{enumerate}
\end{Lemma}
\pf { Part i) is an immediate consequence of the invariance properties of mixed area measures mentioned in Section \ref{sec:pr} (see also
\cite[(5.18)]{Sch1})
and part ii) follows immediately from part i).}
\endpf

\proof[Proof of \rp{decomp}]
Let $K=L+M$ be a decomposable convex body. To see that it is weakly decomposable, notice that by part i) of Lemma~\ref{l:weakly decomposable} 
$$\displaystyle{S_{K+M}=S_{L+2M}=\sum_{r=0}^{n-1}\binom{n-1}{r}2^rS(M[r],L[n-r],\cdot)},$$which is absolutely continuous with respect to 
$$S_K=\sum_{r=0}^{n-1}\binom{n-1}{r}S(M[r],L[n-r],\cdot).$$
\endproof

There are many convex bodies that are not decomposable. For example, every polytope whose 2-dimensional faces are simplices is indecomposable, see \cite[Corollary 3.2.17]{Sch1}. The class of weakly decomposable bodies, however, is much larger than that of decomposable bodies. The next proposition shows that it includes all 
convex polytopes.


\begin{Prop}\label{P:weakly decomposable}
The only weakly indecomposable polytopes are simplices.
\end{Prop}

\begin{pf} In \cite[Lemma 3.1]{SSZ}, we showed that if $K$ is a polytope which is not a simplex, then there exists a convex body $M$ not homothetic to $K$ (obtained from $K$ by moving a facet of $K$ along the direction of its outer normal) such that $S(K[r],M[n-r],\cdot)$ is absolutely continuous with respect to $S_K$, for all $0\leq r<n$. 
Combing this result with Lemma~\ref{l:weakly decomposable}, we obtain the required statement.
\end{pf}

The following theorem, which is the main result of this section, is a generalization of both facts that decomposable convex bodies and polytopes that are not simplices do not satisfy the Bezout inequality (\ref{e:Bezout}) for any $2\leq r\leq n$, see \rt{polytope} and \cite[Theorem 3.3]{SZ}.

\begin{Th}\label{T:weakly}
Let $K\subset \R^n$ be a convex body satisfying 
\begin{equation}\label{e:bezout-wd-r}
V(L_1,\dots,L_r,K[{n-r}])V_n(K)^{r-1}\leq \prod_{i=1}^rV(L_i,K[{n-1}])
\end{equation}
for all convex bodies $L_1,\dots, L_r$ in $\R^n$, where $2\leq r\leq n$. Then $K$ is weakly indecomposable.
\end{Th}

\pf 
It is enough to prove the theorem when $r=2$ as the general case follows from it. Assume that $K$ satisfies 
\begin{equation}\label{e:bezout-wd}
V(L_1,L_2,K[{n-2}])V_n(K)\leq V(L_1,K[{n-1}])V(L_2,K[n-1])
\end{equation}
for all $L_1$, $L_2$ and there exists a convex body~$M$ such that $S_{K+M}$ is absolutely continuous with respect to $S_K$. We need to show that $M$ is homothetic to $K$. The proof combines Theorem \ref{T:dif} and the argument in the proof of \cite[Lemma 3.3]{SSZ}. 

We use induction to show that
\begin{equation}\label{e:A-F-1}
S(M[r],K[n-1-r],\cdot)=\lambda^rS_K
\end{equation}
holds for $0\leq r<n$, where $\lambda=V(M,K[n-1])/V_n(K)$.
The case $r=0$ is trivial. Assume  (\ref{e:A-F-1}) holds for all $0\leq r\leq m$. We need to show that (\ref{e:A-F-1}) holds for $r=m+1$. 
First, we claim that for any convex body $L$ we have
\begin{equation}\label{e:A-F-2}
V(L,M[m+1],K[n-m-2])\leq \lambda^{m+1}V(L,K[n-1]),
\end{equation}
with equality when $L=K$.
If $m=0$, the claim follows immediately from our assumption (\ref{e:bezout-wd}). Therefore we may assume $m\geq 1$. By \re{A-F-1} with $r=m$, for any convex body $L$ we have 
\begin{equation}\label{e:A-F-3}
V(L,M[m],K[n-1-m])=\lambda^mV(L,K[n-1]).
\end{equation}
Applying the Aleksandrov--Fenchel inequality (\ref{eq:af})  to (\ref{e:A-F-3}), we obtain
$$\sqrt{V(L,M[m+1],K[n-2-m])V(L,M[m-1],K[n-m])}\leq \lambda^mV(L,K[n-1]).$$
Then, applying \re{A-F-1} with $r=m-1$ to the second factor in the left-hand side of the above inequality, we obtain (\ref{e:A-F-2}).
Furthermore, if $L=K$, (\ref{e:A-F-3}) produces 
\begin{equation*}
\begin{split}
V(K,M[m+1],K[n-m-2])&=V(M,M[m],K[n-1-m])=\lambda^mV(M,K[n-1])\\
&=\lambda^{m+1}V_n(K)=\lambda^{m+1}V(K,K[n-1]).
\end{split}
\end{equation*}
Thus there is equality in (\ref{e:A-F-2}) when $L=K$.

Next, for a continuous function $f:\mS^{n-1}\to\mathbb{R}$ and sufficiently small $|t|$, define
$$F(t)=V(K_t,M[m+1],K[n-m-2])-\lambda^{m+1}V(K_t,K[n-1]),$$
where $K_t=K_t(f)$ as in \re{perturb}.
By assumption and Lemma~\ref{l:weakly decomposable}, the mixed area measure $S(M[m+1],K[n-r-2],\cdot)$ is absolutely continuous with respect to $S_K$. Therefore, by Theorem \ref{T:dif}, $F(t)$ is differentiable at $t=0$ with
$$F'(0)=\frac{1}{n}\int_{\mS^{n-1}}f(u)dS(M[m+1],K[n-m-2],u)-\lambda^{m+1}\frac{1}{n}\int_{\mS^{n-1}}f(u)dS_K(u).$$
Since $F(0)=0$ and $F(t)\leq 0$ for small $|t|$, it follows that $F'(0)=0$ or equivalently,
$$\int_{\mS^{n-1}}f(u)dS(M[m+1],K[n-m-2],u)=\lambda^{m+1}\int_{\mS^{n-1}}f(u)dS_K(u).$$
But $f$ is arbitrary, which implies that 
$$S(M[m+1],K[n-m-2],\cdot)=\lambda^{m+1}S_K$$
and, hence, (\ref{e:A-F-1}) holds for $r=m+1$. Finally, using (\ref{e:A-F-1}) with $r=n-1$ we see that $S_M$ is proportional to $S_K$. Therefore, by 
Minkowski's Uniqueness Theorem (see \cite[Theorem 8.1.1]{Sch1}) $M$ is homothetic to $K$, as required.
\endpf

 In light of \rp{weakly decomposable} it seems plausible that the only weakly indecomposable convex bodies are simplices.  Then Theorem \ref{T:weakly}  would imply that \rcj{BE} is true. In fact, one can construct many weakly decomposable convex bodies that are not convex polytopes. For example, one can start with a polytope and fix one of its facets. Then any convex body which coincides with the polytope in a neighborhood of the fixed facet will be weakly {decomposable}, as the same argument as in the proof of \rp{weakly decomposable} applies. Still a complete description of {the class of} weakly decomposable bodies is open.

\section{Isomorphic versions of inequalities \re{B} and \re{main}.}\label{S:Iso}
As follows from \rt{n-1}, the inequality \re{B} may not hold when $K$ is not an $n$-simplex. It is natural to ask if \re{B} can be relaxed so that it holds for arbitrary convex sets $L_1,\dots, L_n$ and $K$ in $\R^n$. It turns out that such an inequality is obtained by introducing the constant $n$ in the right-hand side:
 \begin{equation}\label{e:IsoB}
V(L_1,\dots,L_{n})V_n(K)\leq n V(L_1,K[{n-1}])V(L_2,\dots, L_{n},K).
\end{equation}
To show this, we follow the idea of Jian Xiao who used Diskant's inequality to prove that for any convex bodies
$K,L$ in $\R^n$ one has
\begin{equation}\label{e:Xiao}
V_n(K)L\subseteq nV(L,K[n-1])K,
\end{equation}
up to a translation, see \cite[Section 3.1]{Xiao}. Then, \re{IsoB} follows from \re{Xiao} by the monotonicity of the mixed volume.

The inequality \re{IsoB} is, in fact, sharp. For example, one can take $L_1$ to be a unit segment $L_1=[0,u]$ for some $u\in\mS^{n-1}$
and $K$ to be a cylinder $K=L_1\times K'$  for some $(n-1)$-dimensional convex body $K'$ in the orthogonal hyperplane $u^\perp$.
Then $V_n(K)=V_{n-1}(K')=n V(L_1,K[{n-1}])$ by \re{proj}. 
Choose any $(n-1)$-dimensional convex bodies $L_2,\dots, L_n$ in $u^\perp$. 
From the monotonicity of the mixed volume we see that $V(L_2\dots,L_{n},K')=0$ as all of the sets are contained in some $(n-1)$-dimensional ball.
Now, since $K=L_1+K'$, by the linearity and symmetry properties
$$V(L_2,\dots, L_{n},K)=V(L_2\dots,L_{n},L_1)+V(L_2\dots,L_{n},K')=V(L_1,L_2\dots,L_{n}),$$
which provides equality in \re{IsoB}. Note that both sides of the equality are positive. 
{Thus, one can use an approximation argument to show that the constant $n$ in \re{IsoB} cannot be improved for the class of convex bodies as well.}

Next we turn to an isomorphic version of \re{main}. From \re{IsoB} we have
 \begin{equation}\label{e:IsoB-main}
V(L_1,\dots,L_{n-1},K)V_n(K)\leq n V(L_1,K[{n-1}])V(L_2,\dots, L_{n-1},K,K).
\end{equation}
However, we do not expect this inequality to be sharp. Although we do not have a better estimate in general, we can show that for any 
{\it zonoids} $L_1,\dots, L_{n-1}$ and any $K$ in $\R^n$ one has
  \begin{equation}\label{e:IsoB-main-2}
V(L_1,\dots,L_{n-1},K)V_n(K)\leq (n-1) V(L_1,K[{n-1}])V(L_2,\dots, L_{n-1},K,K),
\end{equation}
and in this class the inequality is sharp. (See \cite[p. 191]{Sch1} for the definition of zonoids.) Indeed, similar to the proof of \cite[Theorem 5.6]{SZ}, it is enough to show that \re{IsoB-main-2} holds
when $L_1,\dots, L_{n-1}$ are orthogonal segments. This is a particular case of the Loomis-Whitney type inequalities in \cite[Theorem 1.7]{AAGHV}
(see also \cite{FGM, GHP}). Note that \re{IsoB-main-2} becomes equality, for example, when $L_1=[0,e_1],\dots, L_{n-2}=[0,e_{n-1}]$ 
and $K=\conv\{K',e_n\}$, where $K'$ is the unit $(n-1)$-dimensional cube in $e_n^\perp$.

An isomorphic version of \re{Bezout} was first studied in \cite{SZ}. It was shown that there exists a constant $c_{n,r}>0$ such that 
\begin{equation}\label{e:IsoBezout}
V(L_1,\dots,L_r,K[{n-r}])V_n(K)^{r-1}\leq c_{n,r}\prod_{i=1}^r V(L_i,K[{n-1}])
\end{equation}
holds for arbitrary convex bodies $L_1,\dots, L_r$ and $K$ in $\R^n$, see \cite[Theorem 5.7]{SZ}. 
Since then, several new results on estimating the constant $c_{n,r}$ have been obtained, see \cite{AFO, BGL, Xiao}. Moreover, a generalization of \re{IsoBezout} also appeared in \cite{Xiao}. In particular, \cite[Theorem 1.5]{BGL} provides an isomorphic version of inequality \re{Bezout r=2}
\begin{equation}
V(L_1,L_2,K[{n-2}])V_n(K)\leq 2V(L_1,K[{n-1}])V(L_2,K[{n-1}]).\nonumber
\end{equation}
This inequality becomes equality, for example, when $L_1=[0,e_1]$, $L_2=[0,e_2]$, and $K=\conv\{K',e_3,\dots,e_n\}$, where
$K'$ is the unit square in the span of $\{e_1,e_2\}$.







\end{document}